# On interpretation of complex numbers

Alexander I. Zhbanov

Research Center for Applied Sciences, Academia Sinica, 128, Section 2, Academia Road
Nankang, Taipei 115, Taiwan
Electronic address: azhbanov@gate.sinica.edu.tw; azhbanov@yahoo.com

**Abstract**
We suggest new types and interpretation of complex and hypercomplex numbers for which the commutative, associative, and distributive laws and the norm axioms are trivially satisfied.

## 1. Introduction

Any positive or negative real number can be written down in the form of unsigned vector of size 2 containing separately positive and negative parts of this number.

For example negative real number -2.1 is $\left\{\begin{matrix}0\\2.1\end{matrix}\right\} = \left\{\begin{matrix}1.1\\3.2\end{matrix}\right\} = \left\{\begin{matrix}4\\6.1\end{matrix}\right\} = ...$

Let's call the notation of number -2.1 in the form $\left\{\begin{matrix}0\\2.1\end{matrix}\right\}$ the reduced form.

In a reduced form of an unsigned vector one of its components is equal to zero.

The arbitrary unsigned vector is $\left\{\begin{matrix}a\\b\end{matrix}\right\}$ where $a$ and $b$ are positive real numbers.

Graphical representation of unsigned vector is shown in Fig. 1.

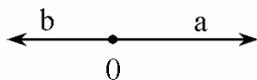

Fig. 1. Unsigned vector

Addition of two real numbers represented in the form of unsigned vectors:
$$\left\{\begin{matrix}a_1\\b_1\end{matrix}\right\} + \left\{\begin{matrix}a_2\\b_2\end{matrix}\right\} = \left\{\begin{matrix}a_1+a_2\\b_1+b_2\end{matrix}\right\}.$$

Multiplication of unsigned vectors:
$$\left\{\begin{matrix}a_1\\b_1\end{matrix}\right\} \times \left\{\begin{matrix}a_2\\b_2\end{matrix}\right\} = \left\{\begin{matrix}a_1 \times a_2 + b_1 \times b_2\\a_1 \times b_2 + a_2 \times b_1\end{matrix}\right\}.$$

Table of multiplication

| × | A | b |
|---|---|---|
| a | a | b |
| b | b | a |

Zeros of unsigned vectors are $\left\{\begin{matrix}0\\0\end{matrix}\right\} = \left\{\begin{matrix}1\\1\end{matrix}\right\} = \left\{\begin{matrix}4\\4\end{matrix}\right\} = ...$

Unit elements of unsigned vectors: $\left\{\begin{matrix}1\\0\end{matrix}\right\} = \left\{\begin{matrix}2\\1\end{matrix}\right\} = \left\{\begin{matrix}4\\3\end{matrix}\right\} = ...$

As an example let's multiply two unsigned vectors +2 and -2 written down in arbitrary form:

$$\begin{Bmatrix}3\\1\end{Bmatrix} \times \begin{Bmatrix}4\\6\end{Bmatrix} = \begin{Bmatrix}3\times 4 + 1\times 6\\3\times 6 + 1\times 4\end{Bmatrix} = \begin{Bmatrix}18\\22\end{Bmatrix} = \begin{Bmatrix}0\\4\end{Bmatrix}.$$

It is obvious that representation of a real number in the form of reduced unsigned vector does not give anything new. The constant component of an unsigned vector contains any insignificant information on its prehistory.

## 2. Complex numbers

Whether three-sign real numbers are possible? Whether are possible an unsigned vectors of dimension 3?

Let's call vector $\begin{Bmatrix}a\\b\\c\end{Bmatrix}$ as the (3)-vector where $a, b, c \geq 0$.

The reduced forms of (3)-vector are $\begin{Bmatrix}a\\0\\0\end{Bmatrix}, \begin{Bmatrix}0\\b\\0\end{Bmatrix}, \begin{Bmatrix}0\\0\\c\end{Bmatrix}, \begin{Bmatrix}a\\b\\0\end{Bmatrix}, \begin{Bmatrix}a\\0\\c\end{Bmatrix}$ or $\begin{Bmatrix}a\\b\\0\end{Bmatrix}$.

For reduction of (3)-vectors to a reduced form it is necessary to decrease all its components by the minimal one.

Example

$$\begin{Bmatrix}2\\3\\1\end{Bmatrix} = \begin{Bmatrix}1\\2\\0\end{Bmatrix}.$$

By obvious way we define a rule of addition of two (3)-vectors:

$$\begin{Bmatrix}a_1\\b_1\\c_1\end{Bmatrix} + \begin{Bmatrix}a_2\\b_2\\c_2\end{Bmatrix} = \begin{Bmatrix}a_1+a_2\\b_1+b_2\\c_1+c_2\end{Bmatrix}$$

Let's define the rule of multiplication of two (3)-vectors:

$$\begin{Bmatrix}a_1\\b_1\\c_1\end{Bmatrix} \times \begin{Bmatrix}a_2\\b_2\\c_2\end{Bmatrix} = \begin{Bmatrix}a_1\times a_2 + b_1\times c_2 + b_2\times c_1\\a_1\times b_2 + a_2\times b_1 + c_1\times c_2\\a_1\times c_2 + a_2\times c_1 + b_1\times b_2\end{Bmatrix}.$$

Table of multiplication

| × | a | B | c |
|---|---|---|---|
| a | a | b | c |
| b | b | c | a |
| c | c | a | b |

Example of multiplication of two (3)-vectors:

$$\begin{Bmatrix}2\\1\\0\end{Bmatrix} \times \begin{Bmatrix}0\\2\\1\end{Bmatrix} = \begin{Bmatrix}2\times 0 + 1\times 1 + 2\times 0\\2\times 2 + 0\times 1 + 0\times 1\\2\times 1 + 0\times 0 + 1\times 2\end{Bmatrix} = \begin{Bmatrix}1\\4\\4\end{Bmatrix} = \begin{Bmatrix}0\\3\\3\end{Bmatrix}.$$

Two examples of raising of (3)-vectors to the second power are the following.

First example:
$$\left\{\begin{matrix}1\\1\\0\end{matrix}\right\}^2 = \left\{\begin{matrix}1\\2\\1\end{matrix}\right\} = \left\{\begin{matrix}0\\1\\0\end{matrix}\right\}.$$

Second example:
$$\left\{\begin{matrix}1\\2\\0\end{matrix}\right\}^2 = \left\{\begin{matrix}1\\4\\4\end{matrix}\right\} = \left\{\begin{matrix}0\\3\\3\end{matrix}\right\}.$$

Graphical representation of (3)-vector is shown in Fig. 2.

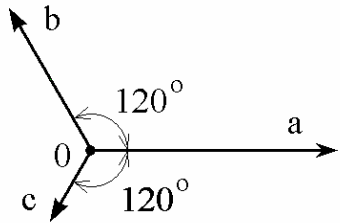

Fig. 2. (3)-vector can be represented as the sum of three vectors on a plane.

Zeros of (3)-vectors are $\left\{\begin{matrix}0\\0\\0\end{matrix}\right\} = \left\{\begin{matrix}1\\1\\1\end{matrix}\right\} = \left\{\begin{matrix}5\\5\\5\end{matrix}\right\} = ...$

It is possible to define the norm of (3)-vector by the law of cosine.
For example:
$$\left\|\left\{\begin{matrix}a\\b\\0\end{matrix}\right\}\right\| = \sqrt{a^2 + b^2 - a \times b}.$$

General case:
$$\left\|\left\{\begin{matrix}a\\b\\c\end{matrix}\right\}\right\| = \sqrt{a^2 + b^2 + c^2 - a \times b - a \times c - b \times c}.$$

It is easy to prove that (3)-vectors and complex numbers are equivalent.
We will not do it. Let's simply show two examples.

Let's raise a complex number $z = x + i \times y = 0 + i \times 1 = i$ to the second power. We receive $z^2 = -1$. We have similar for a (3)-vector in square:
$$\left\{\begin{matrix}1/\sqrt{3}\\2/\sqrt{3}\\0\end{matrix}\right\}^2 = \left\{\begin{matrix}1/3\\4/3\\4/3\end{matrix}\right\} = \left\{\begin{matrix}0\\1\\1\end{matrix}\right\}.$$

Multiplication of two complex conjugate is real number: $(4 + i \times 3) \times (4 - i \times 3) = 25$.
Analogically:

$$\begin{Bmatrix} 4+3/\sqrt{3} \\ 6/\sqrt{3} \\ 0 \end{Bmatrix} \times \begin{Bmatrix} 4+3/\sqrt{3} \\ 0 \\ 6/\sqrt{3} \end{Bmatrix} = \begin{Bmatrix} 16+24/\sqrt{3}+3+12 \\ 6+24/\sqrt{3} \\ 6+24/\sqrt{3} \end{Bmatrix} = \begin{Bmatrix} 25 \\ 0 \\ 0 \end{Bmatrix}.$$

We never met a complex numbers (see for example [1-4]) in such representation as a (3)-vectors. It is difficult to approve that (3)-vectors are more convenient or evident form of representation of complex numbers. Nevertheless its have some elegance.

## 3. (3*3)-matrixes

Now we shall imagine that each positive real number from (3)-vectors itself represents (3)-vector. Thus we define the (3*3)-matrix.

For example we have first (3*3)-matrix $\begin{bmatrix} 0 & 3 & 2 \\ 2 & 2 & 0 \\ 1 & 0 & 3 \end{bmatrix}$ which consists from three (3)-vectors

$$A_1 = \begin{Bmatrix} 0 \\ 2 \\ 1 \end{Bmatrix},\ B_1 = \begin{Bmatrix} 3 \\ 2 \\ 0 \end{Bmatrix},\ \text{and}\ C_1 = \begin{Bmatrix} 2 \\ 0 \\ 3 \end{Bmatrix}.$$

Let we have second (3*3)-matrix $\begin{bmatrix} 0 & 3 & 1 \\ 1 & 0 & 0 \\ 3 & 2 & 1 \end{bmatrix}$ consisting from $A_2 = \begin{Bmatrix} 0 \\ 1 \\ 3 \end{Bmatrix},\ B_2 = \begin{Bmatrix} 3 \\ 0 \\ 2 \end{Bmatrix},\ \text{и}\ C_2 = \begin{Bmatrix} 1 \\ 0 \\ 1 \end{Bmatrix}.$

Multiplication of two (3*3)-matrixes $[A_1\ B_1\ C_1] \times [A_2\ B_2\ C_2] = [A_3\ B_3\ C_3]$ will be as the following:

$A_3 = A_1 \times A_1 + B_1 \times C_2 + B_2 \times C_1$,
$B_3 = A_1 \times B_2 + A_2 \times B_1 + C_1 \times C_2$,
$C_3 = A_1 \times C_2 + A_2 \times C_1 + B_1 \times B_2$.

Let's multiply these (3*3)-matrixes:

$$\begin{bmatrix} 0 & 3 & 2 \\ 2 & 2 & 0 \\ 1 & 0 & 3 \end{bmatrix} \times \begin{bmatrix} 0 & 3 & 1 \\ 1 & 0 & 0 \\ 3 & 2 & 1 \end{bmatrix} = \begin{bmatrix} 18 & 12 & 18 \\ 11 & 14 & 20 \\ 18 & 19 & 13 \end{bmatrix} = \begin{bmatrix} 7 & 0 & 5 \\ 0 & 2 & 7 \\ 7 & 7 & 0 \end{bmatrix}.$$

Definition of norm of (3*3)-matrix:

$$\|[A\ B\ C]\| = \left\| \begin{bmatrix} a_A & a_B & a_C \\ b_A & b_B & b_C \\ c_A & c_B & c_C \end{bmatrix} \right\| \stackrel{def}{=} [(a_A+a_B+a_C)^2 + (b_A+b_B+b_C)^2 + (c_A+c_B+c_C)^2$$

$$-(a_A+a_B+a_C)\times(b_A+b_B+b_C) - (a_A+a_B+a_C)\times(c_A+c_B+c_C) - (b_A+b_B+b_C)\times(c_A+c_B+c_C)]^{1/2}$$

Examples:

$$\left\| \begin{bmatrix} 0 & 3 & 2 \\ 2 & 2 & 0 \\ 1 & 0 & 3 \end{bmatrix} \right\| = \left\| \begin{Bmatrix} 5 \\ 4 \\ 4 \end{Bmatrix} \right\| = \sqrt{5^2 + 4^2 + 4^2 - 5 \times 4 - 5 \times 4 - 4 \times 4} = 1,$$

$$\left\| \begin{bmatrix} 0 & 3 & 1 \\ 1 & 0 & 0 \\ 3 & 2 & 1 \end{bmatrix} \right\| = \left\| \begin{Bmatrix} 4 \\ 1 \\ 6 \end{Bmatrix} \right\| = \left\| \begin{Bmatrix} 3 \\ 0 \\ 5 \end{Bmatrix} \right\| = \sqrt{3^2 + 5^2 - 3 \times 5} = \sqrt{19},$$

$$\left\| \begin{bmatrix} 7 & 0 & 5 \\ 0 & 2 & 7 \\ 7 & 7 & 0 \end{bmatrix} \right\| = \left\| \begin{Bmatrix} 3 \\ 0 \\ 5 \end{Bmatrix} \right\| = \sqrt{19}.$$

The offered (3*3)-matrixes are highly amusing. Its in many respects differ from hypercomplex numbers [5].

Let's list the basic results for (3*3)-matrixes.

We have two sets of zeros for (3*3)-matrixes.

The "absolute" zeros are $\begin{bmatrix} 0 & 0 & 0 \\ 0 & 0 & 0 \\ 0 & 0 & 0 \end{bmatrix}$ and $\begin{bmatrix} \xi & \zeta & \eta \\ \xi & \zeta & \eta \\ \xi & \zeta & \eta \end{bmatrix}$ where $\xi$, $\zeta$, and $\eta$ are any positive real numbers.

Next set is the zeros from family of rotations.
Examples
$$\begin{bmatrix} \xi & 0 & 0 \\ 0 & \xi & 0 \\ 0 & 0 & \xi \end{bmatrix}, \begin{bmatrix} 0 & 0 & \xi \\ \xi & 0 & 0 \\ 0 & \xi & 0 \end{bmatrix}, \begin{bmatrix} \xi & 0 & 0 \\ 0 & \xi & 0 \\ 0 & \xi & 0 \end{bmatrix}, \ldots$$
There are such 27 zeros. Module of "rotation" zeros is equal to zero. Addition of "rotation" zeros to any (3*3)-matrix does not change its module, but "turns" its elements.

Examples of multiplication of nonzero (3*3)-matrixes by matrixes from family of "rotation" zeros:

$$\begin{bmatrix} 0 & 3 & 2 \\ 2 & 2 & 0 \\ 1 & 0 & 3 \end{bmatrix} \times \begin{bmatrix} 0 & 0 & 1 \\ 0 & 1 & 0 \\ 1 & 0 & 0 \end{bmatrix} = \begin{bmatrix} 8 & 5 & 0 \\ 5 & 0 & 8 \\ 0 & 8 & 5 \end{bmatrix},$$

$$\begin{bmatrix} 0 & 3 & 1 \\ 1 & 0 & 0 \\ 3 & 2 & 1 \end{bmatrix} \times \begin{bmatrix} 0 & 0 & 1 \\ 1 & 0 & 0 \\ 0 & 1 & 0 \end{bmatrix} = \begin{bmatrix} 5 & 3 & 0 \\ 0 & 5 & 3 \\ 3 & 0 & 5 \end{bmatrix},$$

$$\begin{bmatrix} 7 & 0 & 5 \\ 0 & 2 & 7 \\ 7 & 7 & 0 \end{bmatrix} \times \begin{bmatrix} 1 & 0 & 0 \\ 0 & 1 & 0 \\ 0 & 1 & 0 \end{bmatrix} = \begin{bmatrix} 9 & 0 & 12 \\ 0 & 9 & 12 \\ 14 & 7 & 0 \end{bmatrix}.$$

It is visible, that norms of final matrixes are equal to zero.
Obviously norms of resulting matrixes are equal to zero.

Other examples of multiplications:

$$\begin{bmatrix} 0 & 2 & 1 \\ 1 & 0 & 0 \\ 2 & 0 & 0 \end{bmatrix} \times \begin{bmatrix} 0 & 0 & 0 \\ 0 & 0 & 0 \\ 0 & 0 & 1 \end{bmatrix} = \begin{bmatrix} 0 & 0 & 1 \\ 0 & 0 & 2 \\ 2 & 1 & 0 \end{bmatrix},$$

$$\begin{bmatrix} 0 & 2 & 1 \\ 1 & 0 & 0 \\ 2 & 0 & 0 \end{bmatrix} \times \begin{bmatrix} 1 & 0 & 0 \\ 0 & 1 & 0 \\ 0 & 0 & 2 \end{bmatrix} = \begin{bmatrix} 0 & 4 & 3 \\ 2 & 0 & 6 \\ 6 & 3 & 0 \end{bmatrix},$$

$$\begin{bmatrix} 0 & 2 & 1 \\ 1 & 0 & 0 \\ 2 & 0 & 0 \end{bmatrix} \times \begin{bmatrix} 1 & 0 & 0 \\ 0 & 1 & 0 \\ 0 & 0 & 3 \end{bmatrix} = \begin{bmatrix} 0 & 4 & 4 \\ 2 & 0 & 8 \\ 8 & 4 & 0 \end{bmatrix}.$$

It is easy to prove that the square root of 9 can be presented by eight various (3*3)-matrixes.

$$\begin{bmatrix} 3 & 0 & 0 \\ 0 & 0 & 0 \\ 0 & 0 & 0 \end{bmatrix}^2 = \begin{bmatrix} 0 & 0 & 0 \\ 3 & 0 & 0 \\ 3 & 0 & 0 \end{bmatrix}^2 = \begin{bmatrix} 1 & 2 & 2 \\ 0 & 2 & 0 \\ 0 & 0 & 2 \end{bmatrix}^2 = \begin{bmatrix} 0 & 0 & 0 \\ 1 & 0 & 2 \\ 1 & 2 & 0 \end{bmatrix}^2 =$$

$$\begin{bmatrix} 1 & 2 & 2 \\ 0 & 0 & 2 \\ 0 & 2 & 0 \end{bmatrix}^2 = \begin{bmatrix} 0 & 0 & 0 \\ 1 & 2 & 0 \\ 1 & 0 & 2 \end{bmatrix}^2 = \begin{bmatrix} 1 & 0 & 0 \\ 0 & 2 & 2 \\ 0 & 2 & 2 \end{bmatrix}^2 = \begin{bmatrix} 0 & 2 & 2 \\ 1 & 0 & 0 \\ 1 & 0 & 0 \end{bmatrix}^2 = \begin{bmatrix} 9 & 0 & 0 \\ 0 & 0 & 0 \\ 0 & 0 & 0 \end{bmatrix}$$

Let's introduce designations of square roots of +1:

$$1 \stackrel{def}{=} \begin{bmatrix} 1 & 0 & 0 \\ 0 & 0 & 0 \\ 0 & 0 & 0 \end{bmatrix}; \quad -1 \stackrel{def}{=} \begin{bmatrix} 0 & 0 & 0 \\ 1 & 0 & 0 \\ 1 & 0 & 0 \end{bmatrix}; \quad J \stackrel{def}{=} \frac{1}{3} \times \begin{bmatrix} 1 & 2 & 2 \\ 0 & 2 & 0 \\ 0 & 0 & 2 \end{bmatrix}; \quad -J \stackrel{def}{=} \frac{1}{3} \times \begin{bmatrix} 0 & 0 & 0 \\ 1 & 0 & 2 \\ 1 & 2 & 0 \end{bmatrix};$$

$$K \stackrel{def}{=} \frac{1}{3} \times \begin{bmatrix} 1 & 2 & 2 \\ 0 & 0 & 2 \\ 0 & 2 & 0 \end{bmatrix}; \quad -K \stackrel{def}{=} \frac{1}{3} \times \begin{bmatrix} 0 & 0 & 0 \\ 1 & 2 & 0 \\ 1 & 0 & 2 \end{bmatrix}; \quad L \stackrel{def}{=} \frac{1}{3} \times \begin{bmatrix} 1 & 0 & 0 \\ 0 & 2 & 2 \\ 0 & 2 & 2 \end{bmatrix}; \quad -L \stackrel{def}{=} \frac{1}{3} \times \begin{bmatrix} 0 & 2 & 2 \\ 1 & 0 & 0 \\ 1 & 0 & 0 \end{bmatrix}.$$

Table of multiplication of different square root of +1:

| ×  | 1  | -1 | J  | -J | K  | -K | L  | -L |
|----|----|----|----|----|----|----|----|----|
| 1  | 1  | -1 | J  | -J | K  | -K | L  | -L |
| -1 | -1 | 1  | -J | J  | -K | K  | -L | L  |
| J  | J  | -J | 1  | -1 | -L | L  | -K | K  |
| -J | -J | J  | -1 | 1  | L  | -L | K  | -K |
| K  | K  | -K | -L | L  | 1  | -1 | -J | J  |
| -K | -K | K  | L  | -L | -1 | 1  | J  | -J |
| L  | L  | -L | -K | K  | -J | J  | 1  | -1 |
| -L | -L | L  | K  | -K | J  | -J | -1 | 1  |

It is easy to make sure that

$$\begin{bmatrix} 3 & 0 & 0 \\ 6 & 0 & 0 \\ 0 & 0 & 0 \end{bmatrix}^2 = \begin{bmatrix} 3 & 0 & 0 \\ 0 & 0 & 0 \\ 6 & 0 & 0 \end{bmatrix}^2 = \begin{bmatrix} 1 & 0 & 4 \\ 2 & 4 & 2 \\ 0 & 2 & 0 \end{bmatrix}^2 = \begin{bmatrix} 1 & 4 & 0 \\ 0 & 0 & 2 \\ 2 & 2 & 4 \end{bmatrix}^2 =$$

$$\begin{bmatrix} 1 & 4 & 0 \\ 2 & 2 & 4 \\ 0 & 0 & 2 \end{bmatrix}^2 = \begin{bmatrix} 1 & 0 & 4 \\ 0 & 2 & 0 \\ 2 & 4 & 2 \end{bmatrix}^2 = \begin{bmatrix} 1 & 2 & 2 \\ 2 & 0 & 0 \\ 0 & 4 & 4 \end{bmatrix}^2 = \begin{bmatrix} 1 & 2 & 2 \\ 0 & 4 & 4 \\ 2 & 0 & 0 \end{bmatrix}^2 = \begin{bmatrix} 0 & 0 & 0 \\ 27 & 0 & 0 \\ 27 & 0 & 0 \end{bmatrix}$$

Thus the imaginary unit can be presented by eight various (3*3)-matrixes.
Let's introduce designations of square roots of -1:

$$i \stackrel{def}{=} \frac{\sqrt{3}}{3} \times \begin{bmatrix} 1 & 0 & 0 \\ 2 & 0 & 0 \\ 0 & 0 & 0 \end{bmatrix}; \; -i \stackrel{def}{=} \frac{\sqrt{3}}{3} \times \begin{bmatrix} 1 & 0 & 0 \\ 0 & 0 & 0 \\ 2 & 0 & 0 \end{bmatrix}; \; j \stackrel{def}{=} \frac{\sqrt{3}}{9} \times \begin{bmatrix} 1 & 0 & 4 \\ 2 & 4 & 2 \\ 0 & 2 & 0 \end{bmatrix}; \; -j \stackrel{def}{=} \frac{\sqrt{3}}{9} \times \begin{bmatrix} 1 & 4 & 0 \\ 0 & 0 & 2 \\ 2 & 2 & 4 \end{bmatrix};$$

$$k \stackrel{def}{=} \frac{\sqrt{3}}{9} \times \begin{bmatrix} 1 & 4 & 0 \\ 2 & 2 & 4 \\ 0 & 0 & 2 \end{bmatrix}; \; -k \stackrel{def}{=} \frac{\sqrt{3}}{9} \times \begin{bmatrix} 1 & 0 & 4 \\ 0 & 2 & 0 \\ 2 & 4 & 2 \end{bmatrix}; \; l \stackrel{def}{=} \frac{\sqrt{3}}{9} \times \begin{bmatrix} 1 & 2 & 2 \\ 2 & 0 & 0 \\ 0 & 4 & 4 \end{bmatrix}; \; -l \stackrel{def}{=} \frac{\sqrt{3}}{9} \times \begin{bmatrix} 1 & 2 & 2 \\ 0 & 4 & 4 \\ 2 & 0 & 0 \end{bmatrix}.$$

Multiplication of different imaginary units:

| × | i | -i | j | -j | k | -k | l | -l |
|---|---|---|---|---|---|---|---|---|
| i | -1 | 1 | -J | J | -K | K | -L | L |
| -i | 1 | -1 | J | -J | K | -K | L | -L |
| j | -J | J | -1 | 1 | L | -L | K | -K |
| -j | J | -J | 1 | -1 | -L | L | -K | K |
| k | -K | K | L | -L | -1 | 1 | J | -J |
| -k | K | -K | -L | L | 1 | -1 | -J | J |
| l | -L | L | K | -K | J | -J | -1 | 1 |
| -l | L | -L | -K | K | -J | J | 1 | -1 |

Multiplication of imaginary units by square roots of 1:

| × | i | -i | j | -j | K | -k | l | -l |
|---|---|---|---|---|---|---|---|---|
| 1 | i | -i | j | -j | K | -k | l | -l |
| -1 | -i | i | -j | j | -k | k | -l | l |
| J | j | -j | i | -i | -l | l | -k | k |
| -J | -j | j | -i | i | L | -l | k | -k |
| K | k | -k | -l | l | i | -i | -j | j |
| -K | -k | k | l | -l | -i | i | j | -j |
| L | l | -l | -k | k | -j | j | i | -i |
| -L | -l | l | k | -k | j | -j | -i | i |

## Conclusion

In this letter we suggest new interpretation of complex numbers in the form of unsigned vectors of size 3. We introduce new types of hypercomplex numbers in the form of unsigned 3*3 matrixes. For such hypercomplex numbers the commutative, associative, and distributive laws and the norm axioms are trivially satisfied. This numbers have 27 "rotation" zeros, 8 image units, and 8 square roots of +1.